# The distribution of the non-prime numbers
# - A new Sieve -

FABIO GIRALDO-FRANCO and PHIL DYKE

1. *Introduction*

The aim of this work is to show, by introducing a new sieve (*©reated by Fabio Giraldo-Franco*), how prime numbers and non-prime numbers are related, showing the independent patterns that exist in natural numbers and how all of them can be linked.

2. *The Sieve*

Having found that non-prime numbers and prime numbers are nicely separated, and that non-prime numbers form infinite independent sequences or patterns, this sieve was created.
The most important thing is that, even having infinitely many patterns, with infinitely many non-prime numbers in them, no number will be ever in more than one pattern. It means, that the same number will not be found twice.

However, it is necessary to create a method (sieve) to link all these sequences to be able to list the non-prime numbers and the prime numbers separately. (To create any of these lists we have to establish an upper bound that leads us to work with a finitely many patterns each with a finite number of elements, but still no number will be in more than one pattern.)

In this next paragraph, let us describe how the sieve works.

Define a list, let's name it LIST (this will eventually become the list of prime numbers), of natural numbers from 2 to n (the upper bound) we remove from LIST every even number, except 2 which is by definition a prime number, by using the following expression:

**2*2 + 2*N** where $N \geq 0$ and * denotes multiplication (N is for all calculations the sequence of natural numbers starting with 0, i.e., N = {0, 1, 2, 3, 4,…}). This expression (first calculation) gives the sequence of numbers from 4 to infinity by (step) 2. Then, after removing that sequence of numbers from LIST, the original list will become:
LIST = {2, 3, 5, 7, 9, 11, 13, 15, 17, 19, 21, 23, 25, 27, 29, 31, 33, 35, 37, 39, 41, 43, 45, 47, 49, 51, 53, 55, 57, 59, …}.

Now we use the same kind of expression used before, i.e., J*K + M*N where N ≥ 0. The values in the first calculation were: J = 2 (2 is the first term in LIST), K = 2 (2 is the first term in LIST) and M = 2 (2 is the first term in LIST). Having used, in the first calculation, the first term in LIST now we select the next one, i.e., 3 (3 is the second term in LIST).

The expression now is: **3*3 + 6*N** where N ≥ 0. This expression (second calculation) gives the sequence of numbers from 9 to infinity by (step) 6. All numbers in this sequence are now removed from LIST. In this case the values are: J = 3 (3 is the second term in LIST), K = 3 (3 is the second term in LIST) and M = 6 (6 is the product of 3x2, i.e., the second term in LIST multiplied by the first term in LIST). Then, having used the second term in LIST, we select the next one, i.e., 5 (5 is the third term in LIST. -Note that the number 4 was removed with the first calculation-).

The expression now is: **5*5 + 30*N** where N ≥ 0. This expression (third calculation) gives the sequence of numbers from 25 to infinity by (step) 30. All numbers in this sequence are now removed from LIST. In this case the values are: J = 5 (5 is the third term in LIST), K = 5 (5 is the third term in LIST) and M = 30 (30 is the product of 5x3x2, i.e., the third term in LIST multiplied by the second term in LIST multiplied by the first term in LIST). In this case it is possible to see that, for the first time, the value of M is bigger than the product of J*K, i.e., M>J*K. Given that, the new expression (fourth calculation) will use the last value of K, i.e., 5, the last value of M, i.e., 30, but for J it will use the next term in LIST, i.e., 7 (note that 6 was removed with the first calculation).

Continuing, the expression now is: **7*5 + 30*N** where N ≥ 0. This expression (fourth calculation) gives the sequence of numbers from 35 to infinity by (step) 30. All numbers in this sequence are now removed from LIST. In this case the values are: J = 7 (7 is the fourth term in LIST), K = 5 (5 is the third term in LIST) and M = 30 (30 is the product of 5x3x2, i.e., the third term in LIST multiplied by the second term in LIST multiplied by the first term in LIST). In this case it is possible to see that the value of M is less than the value of J*K, i.e., M<J*K. Given that, K takes the current value of J, i.e., 7.

For the next calculation, the expression is: **7*7 + 210*N** where N ≥ 0. This expression (fifth calculation) gives the sequence of numbers from 49 to infinity in steps of 210. All numbers in this sequence are now removed from LIST. In this case the values are: J = 7 (7 is the fourth term in LIST), K = 7 (7 is the fourth term in LIST) and M = 210 (210 is the product of 7x5x3x2, i.e., the fourth term in LIST multiplied by the third term in LIST multiplied by the second term in LIST multiplied by the first term in LIST). In this case we see that, once again, the value of M is bigger than the value of J*K, i.e., M>J*K. Given this, the new expression **11*7 + 210*N** uses the same value of K as the earlier calculation, i.e., 7; also the same value of M, i.e., 210, but for J it uses the next term in LIST, i.e., 11. Hence we note that 8 and 10 were removed with the first calculation but that 9 was removed with the second calculation. Hence there is no redundancy and the sieve is efficient.

By repeating this process we can find the entire list of prime numbers up to any number we want (our upper bound) simply given by LIST. The list of prime numbers will appear, obviously, in order after removing the terms generated in each calculation.

The next expressions are:
11*7 + 210*N where N ≥ 0.     M>J*K  (Sixth calculation)
13*7 + 210*N where N ≥ 0.     M>J*K  (Seventh calculation)
17*7 + 210*N where N ≥ 0.     M>J*K  (Eighth calculation)
19*7 + 210*N where N ≥ 0.     M>J*K  (Ninth calculation)
23*7 + 210*N where N ≥ 0.     M>J*K  (Tenth calculation)
29*7 + 210*N where N ≥ 0.     M>J*K  (Eleventh calculation)
31*7 + 210*N where N ≥ 0.     M<J*K. Note that in this calculation (Twelfth calculation) M<J*K, then the next calculation will change the values of J, K and M.
11*11 + 2310*N where N≥0. M>J*K.  Here the value of M is the product of 11x7x5x3x2 (the first fifth terms in LIST multiplied together). This expression (thirteenth calculation) gives the sequence of numbers from 121 to infinity by (step) 2310.
13*11 + 2310*N where N ≥ 0.    M>J*K  (fourteenth calculation)
17*11 + 2310*N where N ≥ 0.    M>J*K  (fifteenth calculation)
19*11 + 2310*N where N ≥ 0.    M>J*K  (sixteenth calculation)
23*11 + 2310*N where N ≥ 0.    M>J*K  (seventeenth calculation)
etc.

We can iterate as many times as we want using the simple rules stated above.

Note that: in some calculation, where J = K (having done the previous calculations) every number in LIST less than J*K is a prime.
The following table shows the first 79 calculations and patterns of non-prime numbers. The second table will give some additional patterns.

The relationship between patterns is evident now: it is necessary (in order to create a technique to list the prime numbers) to calculate the first expression to find the first pattern; then to remove that pattern from the original list; then, the new list gives the values for the second expression and from the second calculation we get the third expression and so on.

The following table shows the first 79 calculations and patterns of non-primes. *In column three, where the condition is greater than the expression (M>J*K) has been right justified. This has not been done when the condition is less than in order that these rarer conditions stand out and can be easily seen.*

| Calculation | Expression | Situation | Sequences to be removed from LIST |
|---|---|---|---|
|  | J*K + M*N | N≥0 | The values of n goes to infinity if wanted |
| 1 | 2*2 + 2*N | M<J*K | Seq (4 to n by 2) |
| 2 | 3*3 + 6*N | M<J*K | Seq (9 to n by 6) |
| 3 | 5*5 + 30*N | M>J*K | Seq (25 to n by 30) |
| 4 | 7*5 + 30*N | M<J*K | Seq (35 to n by 30) |
| 5 | 7*7 + 210*N | M>J*K | Seq (49 to n by 210) |
| 6 | 11*7 + 210*N | M>J*K | Seq (77 to n by 210) |
| 7 | 13*7 + 210*N | M>J*K | Seq (91 to n by 210) |
| 8 | 17*7 + 210*N | M>J*K | Seq (119 to n by 210) |
| 9 | 19*7 + 210*N | M>J*K | Seq (133 to n by 210) |
| 10 | 23*7 + 210*N | M>J*K | Seq (161 to n by 210) |
| 11 | 29*7 + 210*N | M>J*K | Seq (203 to n by 210) |
| 12 | 31*7 + 210*N | M<J*K | Seq (217 to n by 210) |
| 13 | 11*11 + 2310*N | M>J*K | Seq (121 to n by 2310) |
| 14 | 13*11 + 2310*N | M>J*K | Seq (143 to n by 2310) |
| 15 | 17*11 + 2310*N | M>J*K | Seq (187 to n by 2310) |
| 16 | 19*11 + 2310*N | M>J*K | Seq (209 to n by 2310) |
| 17 | 23*11 + 2310*N | M>J*K | Seq (253 to n by 2310) |
| 18 | 29*11 + 2310*N | M>J*K | Seq (319 to n by 2310) |
| 19 | 31*11 + 2310*N | M>J*K | Seq (341 to n by 2310) |
| 20 | 37*11 + 2310*N | M>J*K | Seq (407 to n by 2310) |
| 21 | 41*11 + 2310*N | M>J*K | Seq (451 to n by 2310) |
| 22 | 43*11 + 2310*N | M>J*K | Seq (473 to n by 2310) |
| 23 | 47*11 + 2310*N | M>J*K | Seq (517 to n by 2310) |
| 24 | 53*11 + 2310*N | M>J*K | Seq (583 to n by 2310) |
| 25 | 59*11 + 2310*N | M>J*K | Seq (649 to n by 2310) |
| 26 | 61*11 + 2310*N | M>J*K | Seq (671 to n by 2310) |
| 27 | 67*11 + 2310*N | M>J*K | Seq (737 to n by 2310) |
| 28 | 71*11 + 2310*N | M>J*K | Seq (781 to n by 2310) |
| 29 | 73*11 + 2310*N | M>J*K | Seq (803 to n by 2310) |
| 30 | 79*11 + 2310*N | M>J*K | Seq (869 to n by 2310) |
| 31 | 83*11 + 2310*N | M>J*K | Seq (913 to n by 2310) |
| 32 | 89*11 + 2310*N | M>J*K | Seq (979 to n by 2310) |
| 33 | 97*11 + 2310*N | M>J*K | Seq (1067 to n by 2310) |

| Calculation | Expression | Situation | Sequences to be removed from LIST |
|---|---|---|---|
| 34 | 101*11 + 2310*N | M>J*K | Seq (1111 to n by 2310) |
| 35 | 103*11 + 2310*N | M>J*K | Seq (1133 to n by 2310) |
| 36 | 107*11 + 2310*N | M>J*K | Seq (1177 to n by 2310) |
| 37 | 109*11 + 2310*N | M>J*K | Seq (1199 to n by 2310) |
| 38 | 113*11 + 2310*N | M>J*K | Seq (1243 to n by 2310) |
| 39 | 127*11 + 2310*N | M>J*K | Seq (1397 to n by 2310) |
| 40 | 131*11 + 2310*N | M>J*K | Seq (1441 to n by 2310) |
| 41 | 137*11 + 2310*N | M>J*K | Seq (1507 to n by 2310) |
| 42 | 139*11 + 2310*N | M>J*K | Seq (1529 to n by 2310) |
| 43 | 149*11 + 2310*N | M>J*K | Seq (1639 to n by 2310) |
| 44 | 151*11 + 2310*N | M>J*K | Seq (1661 to n by 2310) |
| 45 | 157*11 + 2310*N | M>J*K | Seq (1727 to n by 2310) |
| 46 | 163*11 + 2310*N | M>J*K | Seq (1793 to n by 2310) |
| 47 | 167*11 + 2310*N | M>J*K | Seq (1837 to n by 2310) |
| 48 | 173*11 + 2310*N | M>J*K | Seq (1903 to n by 2310) |
| 49 | 179*11 + 2310*N | M>J*K | Seq (1969 to n by 2310) |
| 50 | 181*11 + 2310*N | M>J*K | Seq (1991 to n by 2310) |
| 51 | 191*11 + 2310*N | M>J*K | Seq (2101 to n by 2310) |
| 52 | 193*11 + 2310*N | M>J*K | Seq (2123 to n by 2310) |
| 53 | 197*11 + 2310*N | M>J*K | Seq (2167 to n by 2310) |
| 54 | 199*11 + 2310*N | M>J*K | Seq (2189 to n by 2310) |
| 55 | 211*11 + 2310*N | M<J*K | Seq (2321 to n by 2310) |
| 56 | 13*13 + 30030*N | M>J*K | Seq (169 to n by 30030) |
| 57 | 17*13 + 30030*N | M>J*K | Seq (221 to n by 30030) |
| 58 | 19*13 + 30030*N | M>J*K | Seq (247 to n by 30030) |
| 59 | 23*13 + 30030*N | M>J*K | Seq (299 to n by 30030) |
| 60 | 29*13 + 30030*N | M>J*K | Seq (377 to n by 30030) |
| 61 | 31*13 + 30030*N | M>J*K | Seq (403 to n by 30030) |
| 62 | 37*13 + 30030*N | M>J*K | Seq (481 to n by 30030) |
| 63 | 41*13 + 30030*N | M>J*K | Seq (533 to n by 30030) |
| 64 | 43*13 + 30030*N | M>J*K | Seq (559 to n by 30030) |
| 65 | 47*13 + 30030*N | M>J*K | Seq (611 to n by 30030) |
| 66 | 53*13 + 30030*N | M>J*K | Seq (689 to n by 30030) |
| 67 | 59*13 + 30030*N | M>J*K | Seq (767 to n by 30030) |
| 68 | 61*13 + 30030*N | M>J*K | Seq (793 to n by 30030) |
| 69 | 67*13 + 30030*N | M>J*K | Seq (871 to n by 30030) |
| 70 | 71*13 + 30030*N | M>J*K | Seq (923 to n by 30030) |
| 71 | 73*13 + 30030*N | M>J*K | Seq (949 to n by 30030) |
| 72 | 79*13 + 30030*N | M>J*K | Seq (1027 to n by 30030) |
| 73 | 83*13 + 30030*N | M>J*K | Seq (1079 to n by 30030) |

| Calculation | Expression | Situation | Sequences to be removed from LIST |
|---|---|---|---|
| 74 | 89*13 + 30030*N | M>J*K | Seq (1157 to n by 30030) |
| 75 | 97*13 + 30030*N | M>J*K | Seq (1261 to n by 30030) |
| 76 | 101*13 + 30030*N | M>J*K | Seq (1313 to n by 30030) |
| 77 | 103*13 + 30030*N | M>J*K | Seq (1339 to n by 30030) |
| 78 | 107*13 + 30030*N | M>J*K | Seq (1391 to n by 30030) |
| 79 | 109*13 + 30030*N | M>J*K | Seq (1417 to n by 30030) |

This table shows some further calculations and patterns of non-primes.

| Expression | Situation | Seq to be removed from LIST |
|---|---|---|
| J*K + M*N | N≥0 | The values of n goes to infinity |
| 17*17 + 510510*N | M>K | Seq (289 to n by 510510) |
| 19*19 + 9699690*N | M>K | Seq (361 to n by 9699690) |
| 23*23 + 223092870*N | M>K | Seq (529 to n by 223092870) |
| 29*29 + 6469693230*N | M>K | Seq (841 to n by 6469693230) |
| 31*31 + 200560490130*N | M>K | Seq (961 to n by 200560490130) |
| 37*37 + 7420738134810*N | M>K | Seq (1369 to n by 7420738134810) |
| 41*41 + 304250263527210*N | M>K | Seq (1681 to n by 304250263527210) |
| 43*43 + 13082761331670030*N | M>K | Seq (1849 to n by 13082761331670030) |
| 47*47 + 614889782588491410*N | M>K | Seq (2209 to n by 614889782588491410) |
| 53*53 + 32589158477190044730*N | M>K | Seq (2809 to n by 32589158477190044730) |
| 101*101+232862364358497360900063316880507363070*N | M>K | Seq (10201 to n by 232862364358497360900063316880507363070) |

3. *Comparison of Sieves*

There are different patterns of non-prime numbers that give the same result: the lists of non-prime numbers and prime numbers separately. However, they are not independent patterns.

For example, the *sieve of Eratosthenes*[1] produces the same result, but the arithmetic progressions used by that sieve are generated by the multiples of the primes. That is: firstly create a list of natural numbers excluding the number 1, then take the number 2, fix it as a prime number and eliminate all its multiples from that list. Then, take the number 3, fix it as a prime number and eliminate all its multiples from the list. Then, take the number 5, fix it as a prime number and eliminate all its multiples from the list. Then, do the same thing with all the numbers remaining in the list.
(The typical eliminations, used by Eratosthenes's sieve, from a given prime p are: 2p, 3p, 4p, … )

The patterns used by the *sieve of Eratosthenes* are:
General Pattern = {a, a + d, a + 2d, a + 3d, a + 4d, a + 5d, …}
where a = 2p and d = p.

Pattern 1: with a = 4 and d = 2,
{4, 6, 8, 10, 12, 14, 16, 18, 20, 22, 24, 26, 28, 30, …}.
Pattern 2: with a = 6 and d = 3,
{6, 9, 12, 15, 18, 21, 24, 27, 30, 33, 36, …}.
Pattern 3: with a = 10 and d = 5,
{10, 15, 20, 25, 30, 35, …}.

Clearly those patterns are not independent because, for instance, the numbers 6, 10, 12, 15, 18, 20, 24 and 30 are in more than one pattern. It means that there exist infinitely many numbers that are part of infinitely many different patterns. (In this case, having an upper bound, we can say that there exist finitely many numbers that are part of finitely many different patterns.)

Finding the primes, but avoiding the situation that each number was eliminated more than once, Euler produced a new sieve[2] that works as follows:

Firstly, make a list of natural numbers excluding the number 1, then take the number 2, fix it as a prime number, then multiply 2 by itself and by the rest of the numbers of the list and finally eliminate all their products from that list. After that, take the number 3 (the next number, after 2, still in the list), fix it as a prime number, then multiply it by itself and by the remaining numbers in the list and finally eliminate all their products from the list. Then, do the same thing with all the numbers remaining in the list.

The calculations used by the *Euler's sieve* are:
1.) Multiplying 2 by itself and by the rest of the numbers of the list, gives:
{4, 6, 8, 10, 12, 14, 16, 18, 20, 22, 24, 26, 28, 30, …}.
2.) Then, multiplying 3 by itself and by the remaining numbers of the list, gives:
{9, 15, 21, 27, 33, 39, 45, 51, 57, 63, 69, 75, 81, …}.
3.) Then, multiplying 5 by itself and by the remaining numbers of the list, gives:
{25, 35, 55, 65, 85, 95, 115, 125, 145, 155, …}.
4.) Then, multiplying 7 by itself and by the remaining numbers of the list, gives:
{49, 77, 91, 119, 133, 161, 203, 217, …}

We can see that from the third calculation, in the Euler's sieve, the list of the products do not follow any regular pattern or any unique arithmetic progression.

To show more clearly the difference between these sieves, we will make the following comparison:

We know that the list of the first 100 natural numbers consist of 74 non-prime numbers, 25 prime numbers and the number 1.
In the method used by the sieve of Eratosthenes, in order to eliminate the 74 non-prime numbers from the list, 113 non-primes were created. This means that 39 numbers between 1 and 100 were created more than once. In conclusion, this sieve has created, to find the first 25 prime numbers, 52.70% repeated non-prime numbers. This is not efficient.

Using the same sieve, let's see what happens when we find the prime numbers located between 1 and 1000:
We found that the list of the first 1000 natural numbers is composed by: 831 non-prime numbers, 168 prime numbers and the number 1.
The sieve of Eratosthenes creates 1549 non-primes in order to eliminate 831 non-primes. This means that 718 non-prime numbers between 1 and 1000 were created more than once. In conclusion, this sieve has created, to find the first 168 prime numbers, 86.40% repeated non-prime numbers. Even more inefficient.
This inefficiency will increase as we calculate more prime numbers.

The sieve created by Euler solves the problem of repeated numbers, but from our point of view it has a bigger problem: it does not follow any regular pattern.

On the other hand, the Fabio's sieve shows, clearly, how all non-prime numbers can be found exactly once through infinite regular patterns. Having done that we can list immediately the prime numbers.
The list of prime numbers is just the inverse of that list of non-prime numbers formed by independent patterns.

We create lists of non-prime numbers to show the duplications (or not) generated in each of the sieves, and, as it has been stated before, the list of the prime numbers, which is the list we really want, is just the inverse of any of those lists of the non-prime numbers.

To make the comparison absolutely clear we will give what we call an industrial example (Numbers are treated as raw material):

*Example 3.1.*

We have got a machine that, programmed with certain algorithm, changes a list of natural numbers into two lists: one of non-prime numbers and one of prime numbers.

We want to process the first 1000 natural numbers:
After starting our machine we feed it with a list of natural numbers.

Let's try first the sieve of Eratosthenes:

We put a list of natural numbers, excluding the number 1, in the machine, as follows:
LIST A: {2, 3, 4, 5, 6, 7, … , 999, 1000}
The LIST A will lose the non-prime numbers during the process to become a list of prime numbers; it means that LIST A is ever changing. The machine picks out from LIST A the non-prime numbers and creates a new list with them.
After the process we receive the lists of prime numbers and non-prime numbers that the machine produces.

This process generates the following results:
- One list with 168 prime numbers.
- One list with 1549 non-prime numbers.
- Given that we received 718 repeated non-primes, we can say that the process, taking just the non-prime numbers into account, generates 86.40% of industrial waste.
- The machine has to be fed only once (with LIST A).

Let's try now the Euler's sieve:

We put a list, let's say LIST 1, of natural numbers, excluding the number 1, into the machine. The list will lose the non-prime numbers during the process to become a list of prime numbers; it means that the list will change. The machine picks out some non-prime numbers of the LIST 1 and keeps them apart. Then the machine gives us a new list; let's say LIST 2. Then we put the LIST 2 in the machine. The machine picks out some non-prime numbers of the LIST 2 and keeps them apart. Then the machine gives us a new list; let's say LIST 3. Then we put the LIST 3 in the machine. The machine picks out some non-prime numbers of the LIST 3 and keeps them apart. Then the machine gives us a new list; let's say LIST 4. Then we repeat this step as many times as necessary until all prime numbers and non-primes are separated.
After the process we receive the lists of prime numbers and non-prime numbers that the machine produces.

This process generates the following results:
- One list with 168 prime numbers.
- One list with 831 non-prime numbers.
- Given that we did not receive any unexpected non-prime number, we can say that the process did not generate any industrial waste.
- The machine has to be fed 11 times (with LIST 1, LIST 2, LIST 3, LIST 4, LIST 5, LIST 6, LIST 7, LIST 8, LIST 9, LIST 10 and LIST 11).

Let's finally try the Fabio's sieve:

We put a list of natural numbers, excluding the number 1, in the machine, as follows:
LIST A: {2, 3, 4, 5, 6, 7, …, 999, 1000}

The LIST A will lose the non-prime numbers during the process (which in this case is given by the method stated in this article) to become a list of prime numbers; it means that LIST A is ever changing. The machine picks out from LIST A the non-prime numbers and creates a new list with them. After the process we receive the lists of prime numbers and non-prime numbers that the machine produces.

This process generates the following results:
- One list with 168 prime numbers.
- One list with 831 non-prime numbers.
- Given that we did not receive any unexpected non-prime number, we can say that the process did not generate any industrial waste.
- The machine has to be fed only once (with LIST A).

Summarising:

As the size of the initial list increases the industrial waste generated by the sieve of Eratosthenes increases. (The machine would invest more and more and more time producing waste.) However, the machine needs to be fed just once.

With Euler's sieve, as the size of the initial list increases the number of times that the machine needs to be fed increases. (We would invest more and more and more time feeding the machine.) However, this sieve will never produce any industrial waste.

Finally, with Fabio's sieve, regardless of the size of the initial list, there will never be any industrial waste and the machine will be fed just once.

*Here is the Algorithm or pseudo code for the new sieve.*

To find all the prime numbers less than or equal to a given integer n, by using the Fabio's sieve, do:
1. Create a list of consecutive integers from 2 to n: {2, 3, 4, ..., n}.
2. Let "p" equal 2, the first prime number.
3. Create the expression J*K + M*N.
4. Initially, let J = p; K = p; M = p and N = {0, 1, 2, 3, 4, 5, …}.
5. Compare M against the product of J times K and do:
    5.1. If M < J*K, do:
        5.1.a. Calculate J*K + M*N and remove from the list all the resulting numbers less than or equal to n.
        5.1.b. Find the first number remaining on the list after "p" (which is the next prime), do M equals M times "this number", i.e., (M = M*this number), and then replace "p" with "this number" (be careful, not with the new value of M).
        5.1.c. Do J = p and K = p.
        5.1.d. Repeat step 5 until $K^2$ (which is the same as $p^2$) is greater than n.
    5.2. If M > J*K, do:
        5.2.a. Calculate J*K + M*N and remove from the list all the resulting numbers less than or equal to n.
        5.2.b. Find the first number remaining on the list after "J" (which is the next prime) and replace J with "this number".
        5.2.c. If J*K > n go to step 5.1.b.
        5.2.d. Repeat step 5 until $K^2$ is greater than n.
6. All the remaining numbers in the list are prime numbers.

*4. Finding the first 25 prime numbers by using the Fabio's sieve.*

We start stating the initial sequence, S, of the positive integers, from 2 to 100, i.e., S = {2, 3, 4, 5, 6, 7, 8, 9, 10, 11, 12, 13, 14, 15, 16, 17, ..., 99, 100}.

Then, by following the algorithm instructions, we found:
Pattern A: with a = 4 and d = 2,
{4, 6, 8, 10, 12, 14, 16, 18, 20, 22, 24, 26, 28, 30, 32, 34, 36, 38, 40, …}.
Removing the numbers in pattern A, from S, we have:
S = {2, 3, 5, 7, 9, 11, 13, 15, 17, 19, 21, 23, 25, 27, 29, 31, 33, 35, 37, 39, 41, 43, 45, 47, 49, 51, 53, 55, 57, 59, 61, 63, 65, 67, …, 95, 97, 99}.

Then, we found Pattern B: with a = 9 and d = 6,
{9, 15, 21, 27, 33, 39, 45, 51, 57, 63, 69, 75, 81, 87, 93, 99, …}.
Removing the numbers in pattern B, from S, we have:

S = {2, 3, 5, 7, 11, 13, 17, 19, 23, 25, 29, 31, 35, 37, 41, 43, 47, 49, 53, 55, 59, 61, 65, 67, 71, 73, 77, 79, 83, 85, 89, 91, 95, 97}.

Then, we found Pattern C: with a = 25 and d = 30,
{25, 55, 85, 115, 145, 175, 205, 235, 265, 295, 325, 355, …}.
Then, removing the numbers in pattern C, from S, we have:
S = {2, 3, 5, 7, 11, 13, 17, 19, 23, 29, 31, 35, 37, 41, 43, 47, 49, 53, 59, 61, 65, 67, 71, 73, 77, 79, 83, 89, 91, 95, 97}.

Then, we found Pattern D: with a = 35 and d = 30,
{35, 65, 95, 125, 155, 185, 215, 245, 275, 305, 335, 365, …}.
Then, removing the numbers in pattern D, from S, we have:
S = {2, 3, 5, 7, 11, 13, 17, 19, 23, 29, 31, 37, 41, 43, 47, 49, 53, 59, 61, 67, 71, 73, 77, 79, 83, 89, 91, 97}.

Then, we found Pattern E: with a = 49 and d = 210,
{49, 259, 469, 679, 889, 1099, 1309, 1519, 1729, 1939, …}.
Then, removing the numbers in pattern E, from S, we have:
S = {2, 3, 5, 7, 11, 13, 17, 19, 23, 29, 31, 37, 41, 43, 47, 53, 59, 61, 67, 71, 73, 77, 79, 83, 89, 91, 97}.

Then, we found Pattern F: with a = 77 and d = 210,
{77, 287, 497, 707, 917, 1127, 1337, 1547, 1757, 1967, …}.
Then, removing the numbers in pattern F, from S, we have:
S = {2, 3, 5, 7, 11, 13, 17, 19, 23, 29, 31, 37, 41, 43, 47, 53, 59, 61, 67, 71, 73, 79, 83, 89, 91, 97}.

Then, we found Pattern G: with a = 91 and d = 210,
{91, 301, 511, 721, 931, 1141, 1351, 1561, 1771, 1981, …}.
Then, removing the numbers in pattern G, from S, we have:
S = {2, 3, 5, 7, 11, 13, 17, 19, 23, 29, 31, 37, 41, 43, 47, 53, 59, 61, 67, 71, 73, 79, 83, 89, 97}.

This way we have found, given by the sequence S, the first 25 prime numbers.

5. *Conclusion*

We can conclude that the sieve presented here gives us a better way to find the prime numbers than the known sieves. The list of prime numbers was found using only regular and independent patterns.

Hopefully, this new approach to find prime numbers can become a new tool to try to prove some of the conjectures that still are unsolved about primes.